\documentclass[12pt]{article}
\usepackage{amsmath}
\usepackage{amssymb}
\usepackage{euscript}

\normalbaselines
\newcommand{\qed}{\sqcap\kern-8pt\sqcup}

\newenvironment{proof}
      {\par\noindent{\it Proof\/: }\nopagebreak\normalsize}%
				  {\linebreak[2]\hspace*{\fill}$\qed$\ifdim\lastskip<12pt
       \removelastskip \penalty-200  \vskip12pt  \fi}

\font\frten=eufm10 at 12pt
\font\freight=eufm10
\font\frsix=eufm8
\newfam\frfam\textfont\frfam=\frten
\scriptfont\frfam=\freight
\scriptscriptfont\frfam=\frsix

\newcommand{\CC}{\mathbb{C}}
\newcommand{\RR}{\mathbb{R}}
\newcommand{\NN}{\mathbb{N}}
\newcommand{\ZZ}{\mathbb{Z}}
\newcommand{\QQ}{\mathbb{Q}}

\newcommand{\GG}{\mathbb{G}}

\newtheorem{thm}{Theorem}[section]
\newtheorem{prop}[thm]{Proposition}
\newtheorem{cor}[thm]{Corollary}
\newtheorem{lem}[thm]{Lemma}
\newtheorem{defn}[thm]{Definition}

\newtheorem{rema}[thm]{Remark}

\newtheorem{app}[thm]{Example}

\def \O{{\cal O}}

\def \A{{\cal A}}

\def \N{{\bf N}}
\def \H{{\cal H}}
\def \Spr{{\rm Spec}_{r}\,}

\def \Sp{{\rm Spec}\,}
\def \Spun{{\rm X}^{(1)}\,}

\begin{document}

\title{Witt groups and torsion Picard groups of smooth real curves}

\author{Jean-Philippe Monnier\\
       {\small D\'epartement de Mathématiques, Universit\'e d'Angers,}\\
{\small 2, Bd. Lavoisier, 49045 Angers cedex 01, France}\\
{\small e-mail: monnier@tonton.univ-angers.fr}}
\date{}
\maketitle

\section*{Introduction}

The Witt group of smooth real projective curves was first computed by
Knebusch in \cite{Kn}. If the curve is not
complete but still smooth, the Witt group is also studied in
\cite{Kn} but not explicitely calculated. 
However, for some precise examples of smooth affine curves,
we may find explicit calculations (\cite{Knu} and \cite{Ay-Oj}).
In this paper, the Witt group of a general smooth curve is 
explicitely calculated in terms
of topological and geometrical invariants of the curve.
My method is strongly inspired by Sujatha's calculation of
the Witt group of a smooth projective real surface and
uses a comparison theorem between the graded Witt
group and the étale cohomology groups established in \cite{Mo}. 

In the second part of the paper, we are interested in the torsion subgroup
of the Picard group (denoted by $Pic_{tors} (X)$) 
of a smooth geometrically connected
(non complete) curve $X$ over a real closed field $R$. 
Let $C$ be the algebraic closure of $R$ and
$X_C :=X\times_{\Sp R} \,\Sp C$.
We compute 
$Pic_{tors} (X)$ and $Pic_{tors} (X_C)$
using the Kummer exact sequence for étale cohomology.
These calculations depend
on a new invariant $\eta (X)\in\NN$ (resp. $\eta (X_C )$) 
which we introduce in this note.
We study relations between $\eta (X)$, $\eta (X_C )$ and 
the level and Pythagoras number of curves using 
new results of Huisman and Mahé \cite{Hu-Ma}.

The last part is devoted to the study of smooth affine hyperelliptic 
curves. For such curves we calculate the Witt group
and the torsion Picard group determining the invariant $\eta$.

\section{Preliminaries}

Let $k$ be a field. By a variety over $k$ we mean a reduced separated scheme
of finite type over $\Sp k$. A curve over $k$ is a variety of dimension $1$.

\subsection{Witt group and Etale cohomology}

Let $X$ be a smooth connected curve over a real closed field $R$ and
$R(X)$ denote its function field. Let $W(X)$ and $W(R(X))$ denote
respectively the Witt ring of $X$ and $R(X)$. We denote the set
of codimension one points of $X$ by $X^{(1)}$. Since $X$
is smooth, an element $x\in X^{(1)}$ gives a second residue
homomorphism $\partial_x$ defined on $W(R(X))$ with value
in the Witt group $W(k(x))$ of the residue field at $x$.
Thus one obtains an exact sequence \cite{CT-Sa}: 
$$0\rightarrow W(X)
\rightarrow W(R(X))\stackrel{\partial=(\partial_x )}{\rightarrow}
\bigoplus_{x\in\Spun }
W(k(x))$$
Let $I(R(X))$ be the ideal of even rank forms and $I^n (R(X))$
denote its powers for $n\geq 0$ ($I^0 (R(X))=W(R(X))$). 
Recall that $I^n (R(X))$ is additively generated by
the set of $n$-fold Pfister forms, i.e forms isometric to
forms of the type
$$<<a_1 ,\ldots ,a_n>>:=<1,a_1 >\otimes\ldots\otimes<1,a_n>,\,\,a_i \in R(X)^*$$
For $n\geq 0$ we write 
$I^n (X)=I^n (R(X))\cap W(X)$. Thus the previous exact sequence induces
the following exact sequence:
$$0\rightarrow I^n(X)\rightarrow
I^n (R(X))\stackrel{\partial}{\rightarrow}
\bigoplus_{x\in\Spun }I^{n-1}(k(x))\,\,(1)$$
Let $X(R)$ denote the set of $R$-rational points of 
$X$ and $Cont(X(R),\ZZ)$ denote the set of continuous maps
from $X(R)$ into $\ZZ$. Thus $Cont(X(R),\ZZ)\simeq\ZZ^s$
with $s$ denoting the number of semi-algebraic connected
components of $X(R)$.
Since $W(X)$ injects into $W(R(X))$,
the total signature homomorphism $\Lambda:  W(X)\rightarrow Cont(X(R),\ZZ)$
has kernel precisely the torsion subgroup $W_t (X)$.
For $n\geq 0$ we set $I_t ^n (X)=I^n (X)\cap W_t (X)$.

Let $\H^n$ denote the sheaf asssociated to the presheaf 
$U\mapsto H_{et} ^n (U)$
where for any scheme $Y$ over a field of characteristic $\not= 2$
$H_{et} ^n (Y)$ denotes the étale cohomology group
$H_{et} ^n (Y,\ZZ /2)$. Recall that there is an exact sequence
\cite[Th. 4.2]{BO}:
$$0\rightarrow H^0 (X,\H^n)\rightarrow H^n (R(X))
\stackrel{\partial}{\rightarrow} \bigoplus_{x\in\Spun }
H^{n-1}(k(x))\,\,(2)$$
where $H^0 (X,\H^n )$ denotes the group of global sections of the sheaf
$\H^n$. Let $H_t^0 (X,\H^n ):=\{\alpha\in\H^0 (X,\H^n )|\,\alpha\cup (-1)^k =0$
for some $k\}$ where $(-1)$ is the non trivial element in 
$H^1 (R)=R^* / R^{*2}$. 
Using the exact sequences $(1)$ and $(2)$ (see \cite{Mo}),
for every $n\geq 0$ we have a well defined homomorphism
$$e_n :I^n (X)\rightarrow H^0 (X,\H^n )$$
with kernel $I^{n+1} (X)$. We denote by 
$e'_n :I^n (X)/I^{n+1}(X) \rightarrow H^0 (X,\H^n )$ the corresponding
injective map.

The following theorem gives an affirmative answer to a global
version of a question on quadratic forms raised by Milnor.
\begin{thm}
\label{iso}
\cite{Mo}, \cite{Su}\\
Let $X$ be a smooth integral curve
over a real closed field $R$. 
Then $$e'=( e'_n ):\bigoplus_{n\geq 0} I^n (X)/I^{n+1} (X)
\rightarrow\bigoplus_{n\geq 0} H^0 (X,\H^n )$$
gives an isomorphism between the graded Witt group and
the graded unramified cohomology group.
Moreover the restriction to the torsion part
$$e'=(e'_n ):\bigoplus_{n\geq 0} I_t ^n (X)/I_t ^{n+1} (X)
\rightarrow\bigoplus_{n\geq 0} H_t ^0 (X,\H^n )$$
is also an isomorphism.
\end{thm}

\subsection{Complexification of real varieties}

We recall the definition of the level.
\begin{defn}
The level of a commutative ring with units $A$ is the smallest integer $n$ such that $-1$
is a sum on $n$ squares in $A$. If $X$ is a variety
over a real closed field $R$, the level of $X$ is the level
of the $R$-algebra $\O(X)$ where $\O$ is the structure sheaf of $X$.
\end{defn}

Let $X$ be a smooth connected variety over a real closed field
$R$. We always write $C=R(\sqrt{-1} )$ for the algebraic closure of $R$ and 
$X_C :=X\times_{\Sp R} \Sp C$. We denote the canonical
morphism $X_C \rightarrow X$ by $\pi$.
We need the following lemma concerning complexification 
of real varieties.
\begin{lem} 
\label{connexe}
Let $X$ be a smooth connected variety over a real closed field
$R$. If $X_C$ is not connected then $X_C$ is a disjoint union 
$X_C=Y\coprod Y'$ where both restrictions $\pi|_Y :Y\rightarrow X$ and
$\pi|_{Y'} :Y'\rightarrow X$ are isomorphisms. In particular, we have
$X(R)=\emptyset$. Moreover
$X_C$ is not connected if and only if the level of $R(X)$ is $1$
if and only if the level of $X$ is $1$.
\end{lem}

\begin{proof}
The first part of the lemma is \cite[Lem. 1.1]{CT-S}, the statement
follows from the fact that $\pi$ is finite and étale. 
If $X_C$ is not connected,
since $\pi|_Y :Y\rightarrow X$ is an isomorphism, obviously 
$-1\in R(X)^{*2}$. If $-1\in R(X)^{*2}$ then $R(X)(\sqrt{-1})$
is a product of two fields and $X_C$ is not connected.

We prove now the last part of the lemma. Using Kummer exact sequence
(see next section) for the field $R(X)$ and $X$, we have a 
commutative diagramm with exact lines
$$\begin{array}{clclclclclclclclclc}
 0& \rightarrow & \O (X)^* /\O (X)^{*2} & \rightarrow & H^1 (X) & 
\rightarrow & Pic_2 (X) & \rightarrow & 0 &\\
 & & \downarrow  & & \downarrow  & & \downarrow &
 & &\\
 0 & \rightarrow & R(X)^* /R(X)^{*2} &\rightarrow &
H^1 (R(X)) &\rightarrow &  0 &
\rightarrow & 0 &
\end{array}$$
Since $X$ is smooth $H^1 (X)=H^0 (X,\H ^1 )$ and by a previous exact sequence
the map $H^1 (X)\rightarrow H^1 (R(X))$ is injective. By snake lemma
the map $\O (X)^* /\O (X)^{*2} \rightarrow R(X)^* /R(X)^{*2}$
is also injective and the proof is done.
\end{proof}

\begin{rema}
{\rm 
The last equivalence in the previous lemma is not valid in the singular case.
Consider the affine curve $X$ with coordinates ring $\RR [x,y]/y^2 +(1+x^2 )^2$.
Then $-1$ is clearly a square in $\RR (X)$ but it is well known that 
the level of $X$ is $3$.}
\end{rema}

By a well-known theorem of Pfister, if $X$ is a smooth connected 
curve over a real closed field $R$ with $X(R)=\emptyset$ then 
the level of $R(X)$ is 1 or $2$. Moreover, if $X$ is not complete 
and geometrically connected then 
$X$ is affine and
the level of $X$ is $\leq 3$ \cite{Ma}.
One gets the following consequence
\begin{cor}
Let $X$ be a smooth connected curve over a real closed field $R$
with $X(R)=\emptyset$. Then $X_C$ is connected if and
only if the level of $R(X)$ is $2$. If $X$ is not complete
then $X_C$ is connected if and only if the level of
$X$ is $2$ or $3$.
\end{cor}

\section{Witt group of real curves}

\subsection{Etale cohomology groups of real curves}

For any smooth variety $X$ over a field $k$ of characteristic $\not=2$,
there is an exact sequence of étale sheaves:
$$0\rightarrow\mu_n \rightarrow\GG_m\stackrel{n}{\rightarrow}\GG_m
\rightarrow 0$$
where $\GG_m$ is the sheaf of units and $\mu_n$ the sheaf of $n^{th}$
roots of unity ($\mu_2$ is isomorphic to $\ZZ /2$).
From the previous exact sequence
one gets the exact sequences
$$ 0\rightarrow \O (X)^* /\O (X)^{*n} \rightarrow H_{et}^1 (X,\mu_n)
\rightarrow Pic_n(X)\rightarrow 0\,\, (3)$$
$$ 0\rightarrow Pic(X)/n \rightarrow H_{et}^2 (X,\mu_n)
\rightarrow  Br_n (X)\rightarrow 0\,\, (4)$$
where $  Pic_n (X)$
is the $n$-torsion subgroup of the Picard group of $X$
and $ Br_n (X)$ is the $n$-torsion subgroup of the cohomological Brauer
group of $X$. 

Let $X$ be a smooth connected real curve over a real closed field $R$.
The sets $X(R)$ and $X(C)=X_C (C)$ (the set of $C$-rational points)
are semi-algebraic spaces over $R$.
The cohomology of these sets considered here is always the 
semi-algebraic cohomology. For $R=\RR$ and $C=\CC$, it coincides
with classical cohomology.
There is an exact sequence of étale sheaves
$$0\rightarrow\ZZ /2\rightarrow\pi_* (\ZZ /2)\rightarrow\ZZ /2\rightarrow 0$$
This sequence gives rise to the long exact sequence of étale cohomology groups
$$\ldots H^i (X) \stackrel{res}{\rightarrow} H^i (X_C )
\stackrel{cor}{\rightarrow} H^i (X)\stackrel{\cup (-1)}{\rightarrow} H^{i+1} (X)
\ldots\,\,(5)$$
where the boundary maps from $H^i (X)$ to $H^{i+1} (X)$ are cup-products 
by $(-1)\in H^1 (X)$.

\subsubsection{The case of complete real curves}

Let $X$ be a smooth connected curve over a real closed field $R$
and assume $X/R$ is complete. We may easily calculate the 
étale cohomology groups $H^i (X)$, $i\geq 0$, using 
the exact sequences $(3)$, $(4)$ and computations in \cite{CT-S}
of the torsion and the cotorsion of $Pic(X)$.
Let $q:=dim_R H^1 (X,\O_X)$ and $s$ denote
the number of semi-algebraic connected components of $X(R)$. 
Then $Pic(X)_{tors}\simeq (\QQ /\ZZ )^q\oplus
(\ZZ /2 )^{s-1}$ if $X(R)\not=\emptyset$, 
$Pic(X)_{tors}\simeq (\QQ /\ZZ )^q$ if $X(R)=\emptyset$ \cite[Th. 1.6]{CT-S}.
Moreover, $Pic(X)/2\simeq (\ZZ /2)^s$ if $X(R)\not=\emptyset$,
$Pic(X)/2\simeq (\ZZ /2)$ if $X(R)=\emptyset$ \cite[Th. 1.3]{CT-S};
and by a theorem of Witt $Br(X)\simeq (\ZZ /2)^s$.
For $i>3$ $H^i (X)\simeq H^0 (X(R),\ZZ /2)\oplus H^1 (X(R),\ZZ /2)\simeq
(\ZZ /2)^{2s}$ \cite[Th. 2.3.1]{CT-S}. Thus we obtain
\begin{prop}
\label{etalconcomp}
Let $X$ be a smooth geometrically connected curve over a real closed field $R$,
with $X/R$ complete.
Let $g$ denote the genus of $X$. Thus
\begin{description}
\item[({\cal{i}})] 
If $X(R)\not=\emptyset$ then
$H^0 (X)=\ZZ/2$, $H^1 (X)=(\ZZ /2)^{g+s}$, 
$H^i (X)=(\ZZ /2)^{2s}$ for $i\geq 2$.
\item[({\cal{ii}})] If $X(R)=\emptyset$ then
$H^0 (X)=\ZZ/2$, $H^1 (X)=(\ZZ /2)^{g+1}$, $H^2 (X)=(\ZZ /2)$,
$H^i (X)=0$ for $i\geq 3$.
\end{description}
\end{prop}

We deal now with the case $X_C$ not connected. 
\begin{prop}
\label{etalnoconcomp}
Let $X$ be a smooth connected curve over a real closed field $R$
such that $X_C$ is not connected and $X/R$ is complete.
Let $g$ denote the genus of $Y$ ($Y$ is given by Lemma \ref{connexe}). 
Thus
$H^0 (X)=\ZZ/2$, $H^1 (X)=(\ZZ /2)^{2g}$, $H^2 (X)=(\ZZ /2)$,
$H^i (X)=0$ for $i\geq 3$.
\end{prop}

\begin{proof}
By Lemma \ref{connexe} $X_C=Y\coprod Y'$ and the restrictions of 
$\pi$ to $Y$ and $Y'$ are isomorphisms. Thus $X$ is a smooth 
complete connected curve over $C$ of genus $g$. Then 
$dim_R H^1 (X,\O_X)=2g$. The group $\O (X)^* /\O (X)^{*2}$ is trivial 
since $\pi: Y\rightarrow X$ is an isomorphism and
$Y$ is complete over $C$. Using the exact sequence (3)
we get $H^1 (X)=(\ZZ /2)^{2g}$. The other étale cohomology groups
could be deduced from the previous remarks.
We may calculate $H^1 (X)$ in a different way.
The exact sequence $(5)$ gives
$$0\rightarrow H^0 (X)\rightarrow H^0 (X_C )\rightarrow H^0(X)\rightarrow 
H^1 (X)$$
$$\rightarrow H^1 (X_C )\rightarrow H^1 (X)\rightarrow H^2 (X)
\rightarrow H^2 (X_C)\rightarrow H^2 (X)\rightarrow 0$$
We already know $H^0 (X)$, $H^2 (X)$. Let $Pic^0 (Y)$ denote
the kernel of the degree map defined on $Pic (Y)$, then 
we have an exact sequence 
$$0\rightarrow Pic^0 (Y)\rightarrow Pic (Y)\stackrel{deg}{\rightarrow}\ZZ
\rightarrow 0$$
It is well known that $Pic^0 (Y)$ is a divisible group, hence
the previous exact sequence splits since a divisible group is an 
injective $\ZZ -$module.
Moreover $Pic^0 (Y)_{tors}=(\QQ /\ZZ )^{2g}$
and $Br(Y)=0$.
Using $(3)$ and $(4)$ one gets $H^1 (Y)=H^1 (Y')= Pic_2 (Y)=(\ZZ /2 )^{2g}$
and $H^2 (Y)=H^2 (Y')=Pic(Y)/2= \ZZ /2  $. Counting dimensions 
in the previous exact sequence $(5)$  we obtain
$H^1 (X)=(\ZZ /2)^{2g}$.
\end{proof}

\subsubsection{The case of non complete real curves}

Let $X$ be a smooth connected curve over a real closed field.
In this section we assume that $X/R$ is not complete.
By Nagata's embedding theorem \cite{Na} and resolution
of singularities one can realize $X$ as an open and dense subvariety
of a complete smooth variety $\bar{X}$ over $R$. Let 
$Z:=\bar{X}\setminus X$, then $Z$ consists of $r$ real points
$\{P_1 ,\ldots ,P_r \}$ and $c$ complex points $\{Q_1 ,\ldots ,Q_c \}$
with the notation that a closed point $P$ is real (resp. complex)
if the residue field at $P$ is $R$ (resp. C). Let $Z_C$
denote $\bar{X}_C \setminus X_C$, then $Z_C$ consists of
$r+2c$ closed points. Let $s$ denote
the number of semi-algebraic connected components of $X(R)$ and $t$
the number of such components which are complete (if
$R=\RR$ it means compact). Then $X(R)$ is topologically 
a disjoint sum of $t$ circles and $r$ open intervals i.e. $s=t+r$.
If $X_C$ is connected then $g$ will denote the genus of $\bar{X}_C$
(or $\bar{X}$).
If $X_C$ is not connected then $g$ will denote the genus of $\bar{Y}$
(see Lemma \ref{connexe}, $\bar{X}_C =\bar{Y}\coprod \bar{Y'}$).
We will keep this notations all along this paper also in the complete
case i.e. $X=\bar{X}$.

\begin{prop}
\label{etalcon}
Let $X$ be a smooth geometrically connected curve over a real closed field $R$
such that $X/R$ is not complete.
Thus
\begin{description}
\item[({\cal{i}})] 
If $X( R)\not=\emptyset$ then
$H^0 (X)=\ZZ/2$, $H^1 (X)=(\ZZ /2)^{g+c+s}$, 
$H^i (X)=(\ZZ /2)^{s+t}$ for $i\geq 2$.
\item[({\cal{ii}})] If $X( R)=\emptyset$ then
$H^0 (X)=\ZZ/2$, $H^1 (X)=(\ZZ /2)^{g+c}$,
$H^i (X)=0$ for $i\geq 2$.
\end{description}
\end{prop}

\begin{proof}
By \cite[Th. 1.3]{CT-S} $Pic (X)=D(X)\oplus (\ZZ /2 )^t $ with $D(X)$
a divisible subgroup, thus $Pic (X)/2 \simeq (\ZZ /2 )^t$. Using 
Witt Theorem on Brauer groups and (4), we may calculate $H^2 (X)$.
For $i\geq 3$ $H^i (X)\simeq H^0 (X(R),\ZZ /2)\oplus H^1 (X(R),\ZZ /2)\simeq
(\ZZ /2)^{s+t}$.

Hence we are reduced to calculate $H^1 (X)$. We first calculate
$H^1 (X_C )$. There is an exact sequence 
$$0\rightarrow H^1 (\bar{X}_C )\rightarrow H^1 (X_C )\rightarrow H^0 (Z_C)
\rightarrow H^2 ( \bar{X}_C )\rightarrow H^2 (X_C )=0\,\,(6)$$
which is part of the Gysin sequence \cite[p. 244]{Mi}.
We have $H^2 (X_C )=0$ using (4) since $Pic (X_C )$ is a divisible group
and  $Br (X_C )=0$. Since $Pic^0 (\bar{X}_C)$ is a divisible group,
using the split exact sequence 
$$0\rightarrow Pic^0 (\bar{X}_C)\rightarrow Pic (\bar{X}_C)
\stackrel{deg}{\rightarrow}\ZZ
\rightarrow 0$$
and (4), we have $H^2 ( \bar{X}_C )\simeq (\ZZ /2 )$.
Since $H^1 (\bar{X}_C )\simeq (\ZZ /2 )^{2g}$ and $H^0 (Z_C )=(\ZZ /2 )^{r+2c}$,
counting dimensions in $(6)$, 
one obtains $$H^1 (X_C )\simeq(\ZZ /2 )^{2g+2c+r-1}$$

Now, using the exact sequence
$$0\rightarrow H^0 (X)\rightarrow H^0 (X_C )\rightarrow H^0(X)\rightarrow 
H^1 (X)$$
$$\rightarrow H^1 (X_C )\rightarrow H^1 (X)\rightarrow H^2 (X)
\rightarrow H^2 (X_C)= 0$$
one gets $H^1 (X)=(\ZZ /2)^{g+c+s}$ if $X(R)\not=\emptyset$,
and $H^1 (X)=(\ZZ /2)^{g+c}$ if $X(R)=\emptyset$.
\end{proof}

\begin{prop}
\label{etalnocon}
Let $X$ be a smooth connected curve over a real closed field $R$
such that $X_C$ is not connected and $X/R$ is not complete.
Thus
$H^0 (X)=\ZZ/2$, $H^1 (X)=(\ZZ /2)^{2g+c-1}$, 
$H^i (X)=0$ for $i\geq 3$.
\end{prop}

\begin{proof}
For $H^i (X)$, $i\geq 3$, the proof of the previous proposition works.
We have a decomposition $X_C=Y\coprod Y'$ as in Lemma \ref{connexe}.
The closed subset $Z_C$ of $\bar{X}_C$ consists of $2c$ closed points.
Again by Lemma \ref{connexe}, 
$\bar{X}_C=\bar{Y}\coprod \bar{Y'}$ with $\bar{Y},\bar{Y'}$ complete
over $C$. Obviously $\bar{Y}\setminus Y$ consists of $c$ closed points
of $Z_C$. Using the exact sequence $(6)$ for $Y$, $\bar{Y}$ and 
$\bar{Y}\setminus Y$, we obtain 
$H^1 (Y )\simeq(\ZZ /2 )^{2g+c-1}$. Thus
$$H^1 (X_C )\simeq(\ZZ /2 )^{4g+2c-2}$$
Then counting dimensions in $(5)$ 
$$0\rightarrow H^0 (X)\rightarrow H^0 (X_C )\rightarrow H^0(X)\rightarrow 
H^1 (X)\rightarrow H^1 (X_C )\rightarrow H^1 (X)\rightarrow H^2 (X)=0$$
we have
$$H^1 (X )\simeq(\ZZ /2 )^{2g+c-1}$$
The result is compatible with the fact $X$ and $Y$ are isomorphic
via $\pi$ (Lemma \ref{connexe}).
\end{proof}

\subsection{Separation of real connected components}

Let $X$ be a smooth connected curve over a real closed field.
Let $Cont(X(R),\ZZ /2)$ be the set of continuous map
from $X(R)$ into $\ZZ /2$. For every $n$ there is a map
$$h_n :H^0 (X,\H^n )\rightarrow Cont(X(R),\ZZ /2)$$
For $\alpha\in  H^0 (X,\H^n )$, $p\in X(R)$, $h_n (\alpha ) (p)$
is the image of $\alpha$ in $H^n (k(p))\simeq \ZZ /2$.
This map was studied in \cite{CT-Pa}.
It is well known that 
$$H_t ^0 (X,\H^n )=ker(h_n )$$

The following result will be very useful in this note.
\begin{lem}
\label{surj}
Let $X$ be a smooth connected curve over a real closed field.
Then the map 
$$h_1 :H^0 (X,\H^1 )\rightarrow Cont(X(R),\ZZ /2)\simeq (\ZZ /2)^s$$
is surjective
\end{lem}

\begin{proof}
We may assume that $X(R)\not=\emptyset$. Let $C_1 ,\ldots ,C_s$
denote the semi-algebraic connected components of $X(R)$.
By a theorem of Knebush \cite{Kn}, there exist
$q_1 ,\ldots,q_s\in W(X)$ such that the signature of
$q_i$, denoted by $\hat{q_i}:=\Lambda (q_i )$, is $2$ on $C_i$
and $0$ outside. Clearly one gets $q_i \in I(X)$ and $q_i \not\in I^2 (X)$
since its signature is not divisible by $4$.
In \cite{Mo} the author has defined homomorphisms
$$sign_n: I^n (X)/ I^{n+1} (X)\rightarrow Cont(X(R),
2^n \ZZ /2^{n+1}\ZZ )\simeq (\ZZ /2)^s $$ Let $q\in I^n (X)$ and 
$\bar{q}$ denote the class of $q$
in $I^n (X)/ I^{n+1} (X)$, then $sign_n (\bar{q})=(1/2^n )\hat{q}$ $mod\, 2$.
We may prove that
$sign_1$ is surjective using the classes of the $q_i$ in $I(X)/I^2 (X)$. 
Moreover the following
diagramm is commutative
$$\begin{array}{clclclclclc}
 I (X)/ I^{2} (X)& \stackrel{sign_1}{\rightarrow} & (\ZZ /2)^s &\\
 \downarrow e'_1 & &\parallel &\\
  H^0(X,\H^1)& \stackrel{h_1}{\rightarrow} & (\ZZ /2)^s &
\end{array}$$
Since $e'_1$ is an isomorphism, one gets the result.
\end{proof}

\subsection{Some topological remarks}

In this section we will assume that $R=\RR$. 

In the proof of Propositions \ref{etalcon} and \ref{etalnocon}, 
we calculate $H^1 (X_{\CC} )$ with the Gysin sequence in étale cohomology.
The following lemma is a topological verification of this computation.

\begin{lem}
\label{topo}
Let $Y$ be a smooth connected curve over $\CC$. 
Assume $Y$ is not complete. Let 
$Y\hookrightarrow \bar{Y}$ be the smooth completion of $Y$,
$k$ be the number of closed points in $\bar{Y}\setminus Y$
and $g$ denote
the genus of $\bar{Y}$.
Then $H^1 (Y )\simeq (\ZZ /2)^{2g+k-1}$
\end{lem}

\begin{proof}
With the previous notations $\bar{Y}(\CC )$ has a structure of
a compact 2-manifold, more precisely
a sphere with $g$ handles. By a comparison theorem \cite[Th. 3.12, p.117]{Mi},
for any finite abelian group $M$ we have 
$H_{et}^i (Y,M)\approx H^i (Y(\CC ),M )$ and 
$H_{et}^i (\bar{Y},M)\approx H^i (\bar{Y}(\CC ),M )$.
We write $\bar{Y}(\CC )\setminus Y(\CC ):=\{ P_1 ,\ldots ,P_k\}$.
Let $U:=\coprod_{i=1}^k U_i$ where $U_i$ is a small closed ball 
of $\bar{Y}(\CC )$ centered at $P_i$.
Then $Y(\CC ) \cap U\simeq\coprod_{i=1}^k S^1 $.
Using Mayer-Vietoris exact sequence,
$$0\rightarrow H^0 (\bar{Y}(\CC ),\ZZ)\rightarrow H^0 (Y(\CC ),\ZZ )\oplus
H^0 (U,\ZZ )\rightarrow H^0 (Y(\CC )\cap U,\ZZ)\rightarrow$$
$$ H^1 (\bar{Y}(\CC ),\ZZ)\rightarrow H^1 (Y(\CC ),\ZZ )\oplus
H^1 (U,\ZZ )\rightarrow$$ $$H^1 (Y(\CC )\cap U,\ZZ)
\rightarrow H^2 (\bar{Y}(\CC ),\ZZ)\rightarrow H^2 (Y(\CC ),\ZZ )\oplus
H^2 (U,\ZZ )=0$$
which gives
$$0\rightarrow \ZZ\rightarrow\ZZ\oplus \ZZ^{k} \rightarrow 
\ZZ^k \rightarrow \ZZ^{2g}\rightarrow$$
$$H^1 (Y(\CC ),\ZZ )\rightarrow \ZZ^{k}\rightarrow\ZZ\rightarrow 0$$
we find $H^1 (Y (\CC ),\ZZ )\simeq \ZZ^{2g+k-1}$. Universal-coefficient
formula and comparison theorem give the result.
\end{proof}

Let $X$ be a smooth connected curve over $\RR$. Let $G=Gal (\CC / \RR )$
and consider the space $X(\CC )$ equipped with the continuous action 
of $G$. Then the quotient space $X(\CC )/G$ is a $2$-manifold.
Let $\beta:X(\CC )\rightarrow X(\CC )/G$ denote the quotient map.
Let $\A$ be a $G$-sheaf of group on $X(\CC )$, then the 
equivariant cohomology groups $H^i (X(\CC );G,\A )$ 
are defined in \cite[Ch. 5]{Gr}. Let consider
the following spectral sequence converging to the equivariant 
cohomology groups 
$$E_2^{p,q}=H^p (X(\CC )/G,\H^q (G,\A))\Rightarrow 
E^{p+q}=H^{p+q} (X(\CC );G,\A)$$
where $\H^q (G,\A)$ is the sheaf on $X(\CC )/G$
associated to the presheaf
$$U\mapsto H^q (\beta^{-1}(U);G,\A)$$
For the sheaf $\A =\ZZ /2$, we have \cite[1,23; 1-24]{Ni}
$$E_2^{p,0}=H^p (X(\CC )/G, \ZZ /2 );\,E_2^{p,q}=H^p (X(\RR ),\ZZ /2)
\,\,\rm{if}\,q>0$$
and $H^n (X)\simeq H^n (X(\CC );G,\ZZ /2)$.
Then the spectral sequence $E_2^{p,q}$ for $\ZZ /2$
consists of the following
$$\begin{array}{clclclclclc}
 E_2^{0,2}=H^0 (X(\RR ),\ZZ /2)& E_2^{1,2}=H^1 (X(\RR ),\ZZ /2) & 
 E_2^{2,2}=0 &0\\
 E_2^{0,1}=H^0 (X(\RR ),\ZZ /2) & E_2^{1,1}=H^1 (X(\RR ),\ZZ /2) &
 E_2^{2,1}=0 &0\\
 E_2^{0,0}=H^0 (X(\CC )/G, \ZZ /2 )& E_2^{1,0}=H^1 (X(\CC )/G, \ZZ /2 ) & 
 E_2^{2,0}=H^2 (X(\CC )/G, \ZZ /2 ) & 0
\end{array}$$
We have clearly $$E_2^{0,0}=E_{\infty}^{0,0}=E^0$$
Moreover $E_2^{1,0}=E_{\infty}^{1,0}$ and $E_{\infty}^{0,1}
=ker(d_2^{0,1}:E_2^{0,1}\rightarrow E_2^{2,0})$.
Hence the filtration $0\subseteq E_1^1\subseteq E^1$ is given by
$E_1^1 =E_{\infty}^{1,0}=H^1 (X(\CC )/G,\ZZ /2 )$ and
and $E_1 /E_1^1 =E_{\infty}^{0,1}
=ker(d_2^{0,1}:H^0 (X(\RR ),\ZZ /2)\rightarrow H^2 (X(\CC )/G, \ZZ /2 ))$.
We have obviously the following exact sequence which is
the five terms exact sequence of low degree
$$0\rightarrow  H^1 (X(\CC )/G,\ZZ /2 )\rightarrow H^1 (X)
\stackrel{e}{\rightarrow} H^0 (X(\RR ),\ZZ /2)
\stackrel{d_2^{0,1}}{\rightarrow} H^2 (X(\CC )/G, \ZZ /2 )\,\,(6)$$
where $e$ is the edge map. In the Bloch-Ogus spectral sequence
$H^p (X,\H^q )\Rightarrow H^{p+q} (X)$, we have $H^p (X,\H^q )=0$
if $p>q$ ($X$ is smooth). Thus the edge map 
$e': H^1 (X)\rightarrow H^0 (X,\H^1 )$ is an isomorphism.
The following diagramm 
is commutative \cite[Rem. 1.8]{Ni}
$$\begin{array}{clclclclclc}
 H^1 (X)& \simeq & 
 H^0 (X(\CC );G ,\ZZ  /2)\\
 \wr\downarrow e'&  & 
 \downarrow e \\
 H^0 (X,\H^1)& \stackrel{h_1}{\rightarrow} & 
 H^0 (X(\RR ),\ZZ  /2)
\end{array}$$
where $h_1$ is the map defined previously.
Using $(6)$ and the fact that $h_1$ is surjective (Lemma \ref{surj}),
we obtain $$dim_{\ZZ /2} (H^1 (X(\CC )/G,\ZZ /2 ))
=dim_{\ZZ /2}(ker(h_1))=H_t ^0 (X,\H^1 )=dim_{\ZZ /2} (H^1 (X)))-s$$

Since $e$ is surjective in the exact sequence $(6)$,
the differential $d_2^{0,1}:H^0 (X(\RR ),\ZZ /2)
\rightarrow H^2 (X(\CC )/G, \ZZ /2 )$ vanishes.
Consequently $E_2^{2,0}=E_{\infty}^{2,0}=H^2 (X(\CC )/G, \ZZ /2 )$.
Moreover $E_2^{1,1}=E_{\infty}^{1,1}=H^1 (X(\RR ),\ZZ /2 )$
and $E_2^{0,2}=E_{\infty}^{0,2}=H^0 (X(\RR ),\ZZ /2 )$.
We have a filtration $0\subseteq E_2^2\subseteq E_1^2 \subseteq  E^2 =H^2 (X)$
with $E_2^2 =E_{\infty}^{2,0}=H^2 (X(\CC )/G, \ZZ /2 )$,
$E_1^2 /E_2^2 =E_{\infty}^{1,1}=H^1 (X(\RR ),\ZZ /2 )$,
$E^2 /E_1^2 =E_{\infty}^{0,2}=H^0 (X(\RR ),\ZZ /2 )$.
Consequently
we have $$dim_{\ZZ /2} (H^2 (X(\CC )/G,\ZZ /2 ))=
dim_{\ZZ /2} (H^2 (X))-s-t$$

We sum up the previous results.

\begin{prop}
Let $X$ be a smooth complete connected curve over $\RR$.
Thus,
\begin{description}
\item[({\cal{i}})] 
If $X(\RR)\not=\emptyset$ then
$H^0 (X(\CC )/G, \ZZ /2 )\simeq \ZZ /2$,
$H^1 (X(\CC )/G, \ZZ /2 )\simeq (\ZZ /2)^g$,
$H^2 (X(\CC )/G, \ZZ /2 )\simeq 0$.
\item[({\cal{ii}})] If $X(\RR)=\emptyset$ and $X_{\CC}$ is connected then
$H^0 (X(\CC )/G, \ZZ /2 )\simeq \ZZ /2$,
$H^1 (X(\CC )/G, \ZZ /2 )\simeq (\ZZ /2)^{g+1}$,
$H^2 (X(\CC )/G, \ZZ /2 )\simeq \ZZ /2$.
\item[({\cal{ii}})] If $X(\RR)=\emptyset$ and $X_{\CC}$ is not connected then
$H^0 (X(\CC )/G, \ZZ /2 )\simeq \ZZ /2$,
$H^1 (X(\CC )/G, \ZZ /2 )\simeq (\ZZ /2)^{2g}$,
$H^2 (X(\CC )/G, \ZZ /2 )\simeq \ZZ /2$.
\end{description}
\end{prop}

\begin{prop}
Let $X$ be a smooth non complete connected curve over $\RR$.
Thus
\begin{description}
\item[({\cal{i}})] 
If $X(\RR)\not=\emptyset$ then
$H^0 (X(\CC )/G, \ZZ /2 )\simeq \ZZ /2$,
$H^1 (X(\CC )/G, \ZZ /2 )\simeq (\ZZ /2)^{g+c}$,
$H^2 (X(\CC )/G, \ZZ /2 )\simeq 0$.
\item[({\cal{ii}})] If $X(\RR)=\emptyset$ and $X_{\CC}$ is connected then
$H^0 (X(\CC )/G, \ZZ /2 )\simeq \ZZ /2$,
$H^1 (X(\CC )/G, \ZZ /2 )\simeq (\ZZ /2)^{g+c}$,
$H^2 (X(\CC )/G, \ZZ /2 )\simeq 0$.
\item[({\cal{ii}})] If $X(\RR)=\emptyset$ and $X_{\CC}$ is not connected then
$H^0 (X(\CC )/G, \ZZ /2 )\simeq \ZZ /2$,
$H^1 (X(\CC )/G, \ZZ /2 )\simeq (\ZZ /2)^{2g+c-1}$,
$H^2 (X(\CC )/G, \ZZ /2 )\simeq 0$.
\end{description}
\end{prop}

\begin{rema}{\rm In fact one could generalize the 
previous calculation over any real closed
field $R$.
Let $j:X_{et}\rightarrow X_b$ be the morphism of site
introduced in \cite{S}. Then the Leray spectral
sequence for $j_*$
$$H^p (X_b , R^q j_* \ZZ/2 )\Rightarrow H^{p+q} (X)$$
together with a comparison theorem, give the exact sequence $(6)$
over $R$ \cite[20-3-1, p.237]{S}.}
\end{rema}

\subsection{Witt groups of real curves}

\begin{thm}
\label{w1}
Let $X$ be a smooth connected curve over a real closed field $R$.
Let $l$ denote the level of $R(X)$ and $u:=dim_{\ZZ /2 } (H^1 (X))$.
\begin{description}
\item[({\cal{i}})] 
If $X(R)\not=\emptyset$ then
$W(X)\simeq \ZZ^s \oplus (\ZZ /2)^{u-s}$.
\item[({\cal{ii}})] If $X(R)=\emptyset$ and $l=2$ then
$W(X)\simeq \ZZ /4 \oplus (\ZZ /2)^{u-1}$
\item[({\cal{ii}})] If $X(R)=\emptyset$ and $l=1$ then
$W(X)\simeq (\ZZ /2)^{u+1}$.
\end{description}
\end{thm}

\begin{proof}
Assume $X(R)\not=\emptyset$. We have
$W(X)=W_t (X)\oplus \ZZ^s$ since $W_t (X)=I_t (X)$ is the kernel
of the total signature homomorphism $\Lambda$. Since 
$I_t^2 (X)=0$, $I_t (X)$ is a group of exponent $2$.
By Theorem \ref{iso}, $I_t (X)\simeq H_t^0 (X,\H^1 )=ker(h_1)$.
Since $h_1$ is surjective (Lemma \ref{surj})
and $H^1 (X)\simeq H^0 (X,\H^1 )$, the proof is done.

Assume $X(R)=\emptyset$. Thus $W(X)=W_t (X)$ and 
$I^2 (X)=0$. Hence $W(X)$ is a group of exponent $4$
and we have to determine $|W(X)|$ and $|2W(X)|$.
By Theorem \ref{iso}, the rank mod 2 homomorphism
$e_0 ':W(X)/I(X)\rightarrow \ZZ /2$ is an isomorphism
and also the discriminant homomorphism
$e_1 ':I(X)\rightarrow H^0 (X,\H^1 )\simeq H^1 (X)$.
Consequently $|W(X)|=|W(X)/I(X)||I(X)|=2^{u+1}$.
To determine $\mid 2W(X)\mid$ we look at the following 
exact sequence:
$$
 0 \rightarrow N\rightarrow  W(X)/ I (X) 
\stackrel{\otimes <1,1>}{\rightarrow}  2W(X) \rightarrow  0$$
Since $W(R (X))/ I (R (X))=W(X)/ I (X)=\ZZ /2$, the only non-zero element
is the class of $<1>$.

If $N\not= 0$ then $<1,1>=0$ in $2W(X)\subseteq I(X)$.
Since we have injections $I(X)\hookrightarrow I(R(X))$,
$H^1 (X)\hookrightarrow H^1 (R(X))$ and
$H^{1} (R(X))\simeq R(X)^* /R(X)^{*2}$, we get
$e_1 (<<1>>)=-1\in R(X)^{*2}$ i.e $l=1$.
Conversely, if $l=1$ then $N\not=0$. So $N\simeq\ZZ /2$. 
Consequently $2W(X)=0$ and 
$$W(X)\simeq (\ZZ /2 )^{u+1}$$

If $N=0$ i.e $l=2$ then 
$|2W(X)|=|W(X)/I(X)|=2$ and 
$$W(X)\simeq \ZZ /4 \oplus (\ZZ /2)^{u-1}$$
\end{proof}

We finally obtain the explicit 
calculation of the Witt group of a smooth connected
real curve using results of the previous section.

\begin{thm}
\label{w2}
Let $X$ be a smooth complete connected curve 
over a real closed field $R$.
\begin{description}
\item[({\cal{i}})] 
If $X(R)\not=\emptyset$ then
$W(X)\simeq \ZZ^s \oplus (\ZZ /2)^{g}$.
\item[({\cal{ii}})] If $X(R)=\emptyset$ and $X_C$ is connected then
$W(X)\simeq \ZZ /4 \oplus (\ZZ /2)^{g}$
\item[({\cal{ii}})] If $X(R)=\emptyset$ and $X_C$ is not connected then
$W(X)\simeq (\ZZ /2)^{2g+1}$.
\end{description}
\end{thm}

In the case $X_C$ is not connected, we obtain the same result
as \cite[cor. 2.1.7, p. 475]{Knu} concerning the Witt group
of a smooth projective curve over $\CC$.

\begin{thm}
\label{w3}
Let $X$ be a smooth connected curve 
over a real closed field $R$. 
Assume $X/R$ is not complete.
\begin{description}
\item[({\cal{i}})] 
If $X(R)\not=\emptyset$ then
$W(X)\simeq \ZZ^s \oplus (\ZZ /2)^{g+c}$.
\item[({\cal{ii}})] If $X(R)=\emptyset$ and $X_C$ is connected then
$W(X)\simeq \ZZ /4 \oplus (\ZZ /2)^{g+c-1}$
\item[({\cal{ii}})] If $X(R)=\emptyset$ and $X_C$ is not connected then
$W(X)\simeq (\ZZ /2)^{2g+c}$.
\end{description}
\end{thm}

\section{Torsion Picard groups of curves}

\subsection{Torsion Picard groups of real curves}

All along this section $X$ will be a smooth geometrically connected
curve over a real closed field $R$. So $g$ will denote the 
genus of $\bar{X}_C$. We keep the notations of the previous section.
If $X$ is complete, it was shown in \cite{CT-S} using Roitman's 
theorem and a trace argument that
$$Pic_{tors}(X)\simeq (\QQ /\ZZ )^g \oplus (\ZZ /2 )^{s-1}\,\,
{\rm if}\,\, X(R)\not=\emptyset$$
and 
$$Pic_{tors}(X)\simeq (\QQ /\ZZ )^g \,\,{\rm if}\,\, X(R)=\emptyset$$

In this section we will assume that $X$ is not complete.
By \cite[Th. 1.3]{CT-S},
$$Pic (X)\simeq D(X) \oplus (\ZZ /2 )^{t}$$
where $D(X)$ is a divisible group.
Recall that  
$Z:=\bar{X}\setminus X$, and $Z$ consists of $r$ real points
$\{P_1 ,\ldots ,P_r \}$ and $c$ complex points $\{Q_1 ,\ldots ,Q_c \}$.
Let us denote $U_n (X):=\O (X)^* /\O (X)^{*n}$ the group
of units modulo $n$. Let $Jac(\bar{X})$ denote the Jacobian variety
of $\bar{X}$, recall that we have an injective map
$$S:Pic^0 (\bar{X})\rightarrow Jac(\bar{X})(R)$$ 
which is surjective if $\bar{X}(R)\not=\emptyset$.

We will now associate an integer $\eta (X)$ to the curve $X$ as follows.
We denote by $Div(X)$ (resp. $Div(\bar{X})$)
the group of divisors on $X$ (resp. $\bar{X}$)
which is the free abelian group on the closed points of $X$
(resp. $\bar{X}$). We denote by $Div_{rat} (X)$ 
(resp. $Div_{rat}(\bar{X})$) the subgroup of divisors 
rationnally equivalent to $0$
i.e the subgroup of principal divisors.
We have well defined homomorphisms 
$div: R(X)\rightarrow Div_{rat} (X)$ and 
$div: R(\bar{X})\rightarrow Div_{rat} (\bar{X})$.
If $D=\sum_{P} n_P P$ is a divisor,
$supp(D)$ is the set of all points $P$ with $n_p\not=0$.
Let $CH_0 (Z)$ be the group of $0$-cycles on $Z$
modulo rational equivalence, then clearly
$CH_0 (Z)$ is just the free abelian group 
on the closed points of $Z$.
Let $A_0 (Z)$ the free subgroup of $CH_0 (Z)$
which consists of divisors of degree $0$.
Observe that for 
a complex point $Q$, the degree
of the associated divisor $Q $ is 2. In the following 
we will also keep the same notations for
$X_C ,\bar{X}_C ,Z_C:=\bar{X}_C \setminus X_C$.
Let $i:Z\hookrightarrow \bar{X}$ 
and $j:X\hookrightarrow \bar{X}$
be the inclusions.
By \cite[Prop. 1.8, Ch. 1]{Fu}, we have an exact sequence
$$ CH_0 (Z)\stackrel{i_*}{\rightarrow}Pic (\bar{X})
\stackrel{j^*}{\rightarrow} Pic (X)\rightarrow 0$$
Let $B(Z):=Ker( i_* )$ then $B(Z)$ is a free
subgroup of $A_0 (Z)$ of rank $m$. We set
$$\eta (X):=m$$
We may see that $\eta (X)\leq r+c-1$. 
Let $D_1 ,\ldots ,D_{\eta (X)}$ be a basis of $B(Z)$.
Since $i_* (D_j)=0 $ in $Pic^0 (\bar{X} )$,
then $i_* (D_j)=div (f_j )$ with $f_j \in R(\bar{X})=R(X)$
and by the following lemma $f_j \in\O (X)^*$. 
\begin{lem}
\label{inversible}
The following sequence 
$$0\rightarrow \O (X)^* \rightarrow R(X)^* \stackrel{div}{\rightarrow}
Div_{rat} (X)$$
is exact i.e $f\in R(X)^*$ lies in $\O(X)^*$ if and only if
$supp(div(f))\subseteq Z$.
\end{lem}

\begin{proof}
Let $G=Gal (C/R)=\{1 ,\sigma\}$. We have a short exact
sequence $$0\rightarrow \O (X_C )^* \rightarrow C(X)^*
\stackrel{div}{\rightarrow}
Div_{rat} (X_C)\rightarrow 0$$
which induces a long exact sequence in Galois cohomology
$$0\rightarrow (\O (X_C )^* )^G \rightarrow (C(X)^* )^G
\stackrel{div}{\rightarrow}
(Div_{rat} (X_C) )^G \rightarrow H^1 (G ,\O (X_C )^* )\rightarrow\ldots$$
Since $(C(X)^* )^G =R(X)^*$ we get $(\O (X_C )^* )^G =\O (X)^* $.
Moreover, $\pi :X_C \rightarrow X$ induces a flat pull-back (see \cite{Fu})
$\pi^* :Div (X)\rightarrow Div(X_C)$ which is an injection and respects 
rational equivalence. We obtain an injection 
$\pi^* :Div_{rat} (X)\rightarrow (Div_{rat} (X_C))^G$. The statement
follows now easily.
\end{proof}

Let 
$\{ f_1 ,\ldots ,f_{\eta (X)} \}$ be the set 
of units in $\O (X)$ which we create with
$D_1 ,\ldots ,D_{\eta (X)}$. In fact, we may 
complete the previous exact sequence in the following way:
$$ 0\rightarrow R^* =\O (\bar{X} )^* \rightarrow 
\O(X )^*\stackrel{\varphi}{\rightarrow} CH_0 (Z)
\stackrel{i_*}{\rightarrow}Pic (\bar{X})
\stackrel{j^*}{\rightarrow} Pic (X)\rightarrow 0$$
where $\varphi$ is the composition of 
$\O (X )^*\hookrightarrow R(X)=R(\bar{X})$
and $R(\bar{X})\stackrel{div}{\rightarrow} Div(\bar{X})$.

\begin{prop}
\label{unite}
With the previous notations and $n>0$, we have
\begin{description}
\item[({\cal{i}})] 
$U_n (X)\simeq (\ZZ /n)^{\eta (X)}\oplus (\ZZ /2)$
if $n$ is even.
\item[({\cal{ii}})] $U_n (X)\simeq (\ZZ /n)^{\eta (X)}$
if $n$ is odd.
\end{description}
\end{prop}

\begin{proof}
We set $m:=\eta (X)$. We fix a basis 
$D_1 ,\ldots ,D_{m}$ of $B(Z)$ and we get 
the associated $f_i \in\O (X)^*$ for $i=1,\ldots ,m$.
We first claim that any $f\in \O (X)^*$ can be 
written uniquely as a product 
$f=a f_1^{n_1}\ldots f_{m}^{n_m}$ with 
$n_j \in\ZZ$ and $a\in R^*$. We look at $f$ as a rational function
on $\bar{X}$, then $div(f)=D\in Div(\bar{X})$ is an integral combination
of points in $\{P_1 ,\ldots ,P_r ,Q_1 ,\ldots ,Q_c \}$ and 
the degree of $D$ is zero. The divisor
$D$ is in fact in $CH_0 (Z)$ and clearly $i_* (D)=0$.
Hence $D$ can be written uniquely as an integral
combination of $D_1 ,\ldots ,D_m$. Thus we have the claim.

The classes of $f_1 ,\ldots ,f_m$ (and the constant
function $-1$ if $n$ is even) in $U_n (X)$ generate
this group. Since $D_1 ,\ldots ,D_{m}$ 
is a basis of $B(Z)$, we may prove that the order of $f_1 ,\ldots ,f_m$
is exactly $n$ in $U_n (X)$ and $m$ (resp. $m+1$)
is the the minimal number of generators in $U_n (X)$
if $n$ is odd (resp. $n$ is even). 
The statement follows now easily.
\end{proof}

We are now able to calculate the torsion Picard group
of a non complete curve.
\begin{thm}
\label{picard}
Let $X$ be a smooth geometrically connected curve
over a real closed field $R$. Assume $X/R$ is not complete.
Then 
$$Pic_{tors} (X)\simeq (\QQ /\ZZ )^{g+r+c-\eta (X)-1}\oplus (\ZZ /2 )^t$$
with $t=0, r=0$ if $X(R)=\emptyset$.
\end{thm}

\begin{proof}

We have $Pic_{tors} (X)\simeq D(X)_{tors}\oplus (\ZZ /2 )^t$
and for any $n>0$ an exact sequence 
$$ 0\rightarrow U_n (X) \rightarrow H_{et}^1 (X,\mu_n)
\rightarrow Pic_n (X)\rightarrow 0\,\, (3)$$
By \cite[p. 226]{S} $H_{et}^1 (X,\mu_n)\simeq (\ZZ /n )^{g+c+r-1} \oplus
(\ZZ /2 )^{t+1}$, if $n$ is odd one has to drop all summands $\ZZ /2$
(For $n=2$ it gives the result of Theorem \ref{etalcon}).
For $n>1$, we denote by $D_n (X)$ the $n$-torsion subgroup of $D(X)$.
Using Proposition \ref{unite} and the exact sequence $(3)$,
the sequence
$$ 0\rightarrow (\ZZ /n)^{\eta (X)}\oplus (\ZZ /2) \rightarrow 
(\ZZ /n )^{g+c+r-1} \oplus
(\ZZ /2 )^{t+1}\rightarrow D_n (X)\oplus (\ZZ /2 )^t\rightarrow 0$$
is exact if $n$ is even, and 
$$ 0\rightarrow (\ZZ /n)^{\eta (X)}\rightarrow 
(\ZZ /n )^{g+c+r-1} \rightarrow D_n (X) \rightarrow 0$$
is exact 
if $n$ is odd.
Then for any prime number $p>1$ we get
$$D_p (X)\simeq (\ZZ /p )^{g+c+r-\eta (X)-1}$$
By a structure theorem on divisible groups
$D(X)_{tors}$ is a direct sum of some quasicyclic $p$-groups $Z(p^{\infty})$
for some primes $p$. The group $Z(p^{\infty})$ could be seen as the $p$-primary
component of $\QQ /\ZZ$. The result on the $p$-torsion part of $D(X)$
implies that, for any prime $p$, we have exactly $g+c+r-\eta (X)-1$ copies
of $Z(p^{\infty})$ in the decomposition of $D(X)_{tors}$ as a direct sum
of quasicyclic groups. Thus the proof is done.
\end{proof}

\subsection{Torsion Picard group of complex curves}

All along this section $R$ will be a real closed field and 
$C=R(\sqrt{-1})$.

Let $Y$ be a smooth connected
curve over $C$. 
If $Y$ is complete then $Pic_{tors} (X)\simeq (\QQ /\ZZ )^{2g}$,
so we will assume that $Y$ is not complete. Let $\bar{Y}$ denote
the smooth completion of $Y$ and assume that $\bar{Y}\setminus Y$
consists on $k$ closed points. In this section, using Lemma \ref{connexe},
we deal with the case of smooth real curves of level 1.
We keep the notations we used
for real curves, in particular for $\eta (Y)$, $U_n (Y)$ and 
$g$ will denote the 
genus of $\bar{Y}$. 

\begin{prop} 
\label{piccomplexe}
For $Y$ satisfying the previous conditions,
$$Pic_{tors} (Y)\simeq (\QQ /\ZZ )^{2g+k-\eta (Y)-1}$$
\end{prop}

\begin{proof}
By a topological computation, we may prove that 
$H^1 (Y(C),\ZZ)=\ZZ^{2g+k-1}$ and then
$H^1 (Y(C),\ZZ /n)=(\ZZ/n)^{2g+k-1}$ for any $n>0$.
By a comparison Theorem of Huber \cite{Hub}
and since $C$ is algebraically closed,
$$H_{et}^1(Y,\mu_n )=(\ZZ/n)^{2g+k-1}$$
for any $n>0$. Arguing as for real curves, for any $n>0$,
$U_n (Y) \simeq (\ZZ /n )^{\eta(Y)}$, thus 
we get an exact sequence 
$$0\rightarrow (\ZZ /n )^{\eta(Y)}\rightarrow (\ZZ/n)^{2g+k-1}
\rightarrow Pic_n (Y)\rightarrow 0$$
If $B$ is an abelian group, $Hom_{\ZZ} (B,\QQ /\ZZ )$ is the Pontrjagin dual
of $B$. Suppose that $B$ is finite, then $B\simeq Hom_{\ZZ} (B,\QQ /\ZZ )$.
According to the above remark, we get another exact sequence
$$0\rightarrow Pic_n (Y)\rightarrow (\ZZ/n)^{2g+k-1}
\rightarrow (\ZZ /n )^{\eta(Y)}\rightarrow 0$$
which is a split exact sequence of $\ZZ /n$-modules.
Hence
$$Pic_n (Y)\simeq (\ZZ /n)^{2g+k-\eta (Y)-1}$$
for any $n>0$. Since $Pic (Y)$ is a divisible group,
using a structure theorem on divisible groups, the statement follows.
\end{proof}

\subsection{Relations between $\eta (X)$, $\eta (X_C )$ and the level 
of curves}

Let $X$ be a smooth geometrically 
connected curve over a real closed field $R$. We assume
that $X$ is not complete.
We denote by $Z_C =\bar{X}_C \setminus X_C$, $Z=\bar{X} \setminus X$,
$B(Z)$ the free group $Ker(CH_0 (Z )\rightarrow
Pic(\bar{X}))$, $B(Z_C)$ the free group $Ker(CH_0 (Z_C )\rightarrow
Pic(\bar{X}_C ))$. We recall that $\eta (X)$ (resp. $\eta (X_C )$)
is the rank of $B(Z)$ (resp. $B(Z_C )$).

We recall that have an injection 
$\pi^* :Div(\bar{X})\hookrightarrow Div(\bar{X}_C)$
(see the proof of Lemma \ref{inversible}):
For $D\in Div(\bar{X})$ we replace a complex point 
$Q\in Supp(D)$ by the sum $Q' +\bar{Q'}$
where $Q' ,\bar{Q' }$ are the two complex points of $\bar{X}_C$
lying over $Q$. 

\begin{prop}
\label{comparaison}
We always have $\eta (X_C )\geq \eta (X)$, moreover
$0\leq \eta (X)\leq c+r-1$ and $0\leq \eta (X_C )\leq r+2c -1$. 
If $c=0$ (i.e. we have only
real points at ``infinity") then $\eta (X)=\eta (X_C)$.
\end{prop}

\begin{proof}
The first assertion is clear. If $c=0$ then $r\geq 1$
since $X$ is not complete.
Since  $Pic (\bar{X})$ injects into $Pic (\bar{X}_C )^{G}$,
the class of a divisor $D\in Div (\bar{X})\subseteq Div (\bar{X}_C )$
is zero in $Pic (\bar{X})$ if only if its class is zero in $Pic (\bar{X}_C)$.
It is the case, in particular, if $supp(D)$ is contained in
$\{ P_1 ,\ldots ,P_r \}=\Bar{X}\setminus X=\Bar{X}_C \setminus X_C $,
which gives the proof.
\end{proof}

\begin{rema}
{\rm If $r+c=1$ then trivially $\eta (X)=0$.}
\end{rema}

Let $G=Gal (C /R )=\{1 ,\sigma \}$. We denote by
$B(Z_C )^{-}$ (resp. $B(Z')^{+}$) the subgroup of $B(Z_C )$ which elements
are anti-invariant (resp. invariant) by $\sigma$
i.e. $B(Z_C )^{-}=Ker( B(Z_C )\stackrel{1+\sigma}{\rightarrow}B(Z_C ))$
(resp. $B(Z_C )^{+}=Ker( B(Z_C )\stackrel{1-\sigma}{\rightarrow}B(Z_C))$). 
Then $B(Z_C )^{-}$ (resp. $B(Z_C )^{+}$) is a free group and we denote 
by $\eta^{-} (X_C )$ (resp. $\eta^{+} (X_C )$) its rank.
Clearly one gets $B(Z_C )^{+}=B(Z)$ and $\eta^{+} (X_C )=\eta (X)$.
Moreover, $B(Z_C )$ could be seen as a subgroup of $Div_{rat} (\bar{X}_C )$
and $B(Z_C )^{-}$ is a subgroup of $Div_{rat} (\bar{X}_C )^{-}=
Ker( Div_{rat} (\bar{X}_C )\stackrel{1+\sigma}{\rightarrow}
Div_{rat} (\bar{X}_C ))$.
Then we obtain a well defined map
$$\phi:  B(Z_C )^{-} \rightarrow H^1 (G, Div_{rat} (\bar{X}_C ))$$
which to $D$ associates the class of $D\in Div_{rat} (\bar{X}_C )$
in the Galois cohomology group $H^1 (G, Div_{rat} (\bar{X}_C ))$.

If $X(R)=\emptyset$, the level of $X$ is $2$ or $3$. 
We reformulate the result of \cite{Hu-Ma} 
in terms of our new invariants.
\begin{prop} 
\label{niv2m}
Assume $X(R)=\emptyset$. If the level of $X$ is $2$
then $\eta (X_C)\geq \eta^{-} (X_C)>0$ and $\eta (X_C)>\eta (X)$.
More precisely the level of $X$ is $2$ if and only if
the map $\phi$ is non zero.
\end{prop}

\begin{proof}
The level of $X$ is $2$ if and 
only there exists $f\in  \O (X_C )^{*}$ such that 
$N(f)=f \bar{f}=-1$.
Assume such $f$ exists and let $D:=Div(f)\in Div(\bar{X}_C )$.
Clearly $D\in B(Z_C)$ and since $D +\bar{D}=0$
one gets $D\in B(Z_C )^{-}$. This proves that $\eta^{-} (X_C )>0$.
 
Conversely assume $D$ lies in $B(Z_C )^{-}$.
It corresponds to the divisor of 
a non trivial element $f_D$ in $\O (X_C )^{*}\subseteq C(X)$.
Moreover $div(f_D \bar{f_D})=D+\bar{D}=0$, so 
$N(f_D )\in R^*$ and we may assume that it is $1$ or $-1$.
We have shown that there is a one-to-one mapping between
$B(Z_C )^{-}$ and the set of functions $f\in  \O (X_C )^{*}$ such that 
$N(f)=f \bar{f}=\pm 1$ modulo $\{ z\in C |z\bar{z} =1\}$.
Hence the level of $X$ is $2$
if and only if there exists $D\in  B(Z_C )^{-}$
such that $N(f_D)=-1$. Since $H^1 (G,C(X)^* )=0$ by Hilbert's Theorem 90,
$N(f_D)=1$ if only if $f_D =g/\bar{g}$ for 
$g\in C(X)^*$ if only if $\phi (D)=0$.

We assume now that the level of $X$ is $2$ and $\eta (X)=\eta (X_C)$.
Since $\eta (X)=\eta^+ (X_C )$, it means that every 
$D\in B(Z_C )^-$ is invariant by $\sigma$ since it can be written
as an integral combination of elements in $B(Z_C )^+$. But then 
$B(Z_C )^-$ is trivial and $\eta^- (X_C )=0$, contradiction.
\end{proof}

\begin{rema}
{\rm 
If $X$ has only one complex point at infinity and $X(R)=\emptyset$, 
such curves are called maximal
in \cite{Hu-Ma}. Then $Z_C =\{Q,\bar{Q} \}$ and
$\eta (X_C)=\eta^- (X_C)=1$ if and only the class of $n(Q -\bar{Q})$
is zero in $Pic^0 (\bar{X_C })$ for a $n\in\NN\setminus\{ 0\}$. We could see 
in \cite{Hu-Ma} 
that there exists maximal curves
with $\eta (X_C)=1$ but the level of $X$ is not $2$.}
\end{rema}

\subsection{Pythagoras number of curves}

We always assume that $X$ is not complete.
For a commutative ring with units $A$, we denote
by $\sum A^2$ (resp. $\sum_{i=1}^n A^2$) the set of sums of 
squares (resp. $n$ squares) in $A$,
$p(A)=\inf\{n\in\N \mid\,\,\sum A^2 =\sum_{i=1}^n A^2 \}$
or $\infty$ the classical Pythagoras number of $A$ and 
$p_* (A)=\inf\{n\in\N \mid\,\,A^*\cap \sum A^2 =\sum_{i=1}^n A^2 \}$
or $\infty$.

Assume $X(R)=\emptyset$.
Let $l:=level (X)$. Since $l$ is finite and 
every element in $\O (X)$ can be written as a difference of two squares
and is finally a sum of squares,
we get:
$$2\leq l\leq p_* (\O (X) )\leq p (\O (X) ) \leq l+1 \leq4$$
Now we would like to know some conditions for which
$$l=p_* (\O (X) )=2$$
Remark that for the field $R(X)$ the equality holds by a result of Pfister.

\begin{prop}
\label{pyth1}
Assume $X(R)=\emptyset$.
We have $l=p_* (\O (X) )=2$
if and only if the map 
$$\pi^* :Pic(X)\rightarrow (Pic (X_C ))^G $$
is surjective.
\end{prop}

\begin{proof}
Consider the short exact sequence
$$0\rightarrow \O (X_C )^* \rightarrow C(X)^*
\stackrel{div}{\rightarrow}
Div_{rat} (X_C)\rightarrow 0$$
which induces a long exact sequence
$$\ldots H^1 (G, C(X)^* )\rightarrow H^1 (G,Div_{rat} (X_C) )
\rightarrow H^2 (G, \O (X_C )^* )\rightarrow H^2 (G, C(X)^* )$$
Now $H^1 (G, C(X)^* )=0$ by Hilbert's Theorem 90 and one gets
an exact sequence
$$0\rightarrow H^1 (G,Div_{rat} (X_C) )\stackrel{\alpha}{\rightarrow}
 \O (X )^* /N(\O (X_C )^* ) \stackrel{\beta}{\rightarrow} R (X )^* /N(C(X)^* )$$
with $N$ denoting the norm.
We have to understand the boundary map $\alpha$.
An element of  $H^1 (G,Div_{rat} (X_C) )$ could be seen as
the class of $D=div(f_D)$ (denoted by  $[D]$) 
with $f_D\in C(X)^*$ and $\bar{D} =-D $.
Then $\alpha ([D] )$ corresponds to the class of $N(f_D)\in \O(X)^*$
modulo $N( \O (X_C )^* )$. We could see that $\alpha$ is well defined.

By a well known Theorem of Pfister, the Pfister form $<<1>>$
is universal over $R(X)$, hence $R (X )^* /N(C(X)^* )=0$
and $\alpha$ is an isomorphism. 

We claim now that $H^1 (G,Div_{rat} (X_C) )$ is the cokernel
of $\pi^* :Pic(X)\rightarrow (Pic (X_C ))^G $ (we don't
have to assume that $X(R)=\emptyset$), and the proof 
will be done. 
The short exact sequence of $G$-modules
$$0\rightarrow Div_{rat} (X_C)\rightarrow Div (X_C)\rightarrow Pic (X_C )
\rightarrow 0$$
gives a long exact sequence 
$$\ldots\rightarrow Div (X)\rightarrow (Pic (X_C ) )^G
\rightarrow H^1 (G,Div_{rat} (X_C) )\rightarrow H^1 (G,Div (X_C) )
\rightarrow\ldots$$
Since $H^1 (G,Div (X_C) )=0$, the statement follows easily.
\end{proof}

\begin{app}
{\rm Let $X$ be the affine plane curve given by the equation 
$x^2 +y^2 +1 =0$. Then $p_* (\O (X) )=2$
since $p(\O (X))=2$ \cite[Th. 3.7]{CDLR}.
Therefore $\pi^* $ is surjective.}
\end{app}

\begin{prop}
\label{pyth2}
We assume $X(R)\not=\emptyset$.
We have $p_* (\O (X) )=2$
if  
$\pi^* :Pic(X)\rightarrow (Pic (X_C ))^G $
is surjective.
Conversely $\pi^*$ is surjective
if $p_* (\O (X) )=2$ and every $f\in \O(X)^*$ 
which is positive on $X(R)$ is a sum of squares in $\O (X)$.
\end{prop}

\begin{proof}
The proof is straightforward 
using the exact sequence
$$0\rightarrow H^1 (G,Div_{rat} (X_C) )\stackrel{\alpha}{\rightarrow}
 \O (X )^* /N(\O (X_C )^* ) \stackrel{\beta}{\rightarrow} R (X )^* /N(C(X)^* )$$
since every positive function 
$f\in R(X)$ on $\Spr R(X)$ is a sum of two squares in $R(X)$.
\end{proof}

\begin{app}
{\rm  
Let $X$ be the affine plane curve given by the
equation $x^2 +y^2 -1 =0$. Since $p(\O (X))=2$ \cite[Th. 3.7]{CDLR},
it follows that $p_* (\O (X) )=2$.
Moreover every positive function in $\O (X)$ is a sum of squares
\cite[Prop. 2.17]{S2}, hence $\pi^*$ is surjective.}
\end{app}

\section{Witt groups and torsion Picard groups of smooth affine plane curves}

\subsection{Conics}

In order to have an easy application of the previous results, we calculate
the Witt group and the torsion Picard group of a smooth curve
in the real plane
given by the zero set of $P(x,y)\in\RR [x,y]$ with $degree(P)\leq 2$.
In any case we have $g=0$.

Up to an isomorphim over $\RR$, 
we are reduced to deal with the following cases:

1) The case of an ellipse: $P(x,y)=x^2 +y^2 -1$.
Then $r=0$, $c=1$, $\eta (X)=0$, $s=t=1$. We have $W(X)\simeq \ZZ\oplus \ZZ /2$
and $Pic_{tors} (X)\simeq \ZZ /2 $ (compare with \cite{Ay-Oj}).

2) The case of a parabola: $P(x,y)=x^2 +y$.
Then $r=1$, $c=0$, $\eta (X)=0$, $s=1$, $t=0$. We have $W(X)\simeq \ZZ$
and $Pic_{tors} (X)\simeq 0$.

3) The case of an hyperbola: $P(x,y)=x^2 -y^2 -1$.
Then $r=2$, $c=0$, $s=2$, $t=0$, $\eta (X)=1$ since the function
$f=(x-y)\in \O (X)^*$ and the divisor of the corresponding rational function
$f=(X-Y)/Z\in \RR (\bar{X})$ has a support contained in the set
of the two infinity points where $\bar{X}$ is the smooth
projective curve given by the equation $X^2 -Y^2 -Z^2 =0$. We have 
$W(X)\simeq \ZZ^2 $
and $Pic_{tors} (X)\simeq 0$ (compare with \cite{Knu}).

4) The case of an imaginary ellipse: $P(x,y)=x^2 +y^2 +1$.
Then $r=0$, $c=1$, $\eta (X)=0$, $s=t=0$. We have $W(X)\simeq \ZZ /4$
and $Pic_{tors} (X)\simeq 0$.

5) The case of a line: $P(x,y)=x$.
Then $r=1$, $c=0$, $\eta (X)=0$, $s=1$. We have $W(X)\simeq \ZZ$
and $Pic_{tors} (X)\simeq 0$

6) The non geometrically connected case: $P(x,y)=x^2 +1$.
Then $r=0$, $c=1$ after a blowing-up at infinity 
of the singular projective curve $X^2 +Z^2 =0$, $s=t=0$.
We have 
$W(X)\simeq \ZZ /2 $
and $Pic_{tors} (X)\simeq 0$.

\subsection{Hyperelliptic curves}

We study here smooth geometrically connected affine curves in 
the real plane given by an equation 
$y^2 +P(x)=0$ with $P(x)\in\RR [X]$ non constant and square free 
($X$ is smooth).
Up to an isomorphism we are reduced to curves with equations
$y^2 +P(x)=0$ and $y^2 -P(x)=0$ with $P(x)$ monic.

\subsubsection{Curves with equations $y^2 +P(x)=0$, $P$ monic}

Let $X$ be the plane curve with equation $y^2 +P(x)=0$. 
We denote by $d$ the degree of $P$ and $k$ the number of real roots 
of $P$. 

First assume that $d$ is even, $d:=2d'$. After some blowings-up at infinity
of the singular projective curve associated to $X$, we see that 
$\bar{X}\setminus X$ consists on one complex point, it means that
$r=0$, $c=1$, $\eta (X)=0$. Moreover by Hurwitz formula $g=d'-1$. The number 
of real roots of $P$ is even ($k=2k'$) and $s=t=k'$. 
Then $$W(X)\simeq \ZZ^{k'}\oplus (\ZZ /2 )^{d'}\,\,{\rm if}\,\, k>0$$
$$W(X)\simeq \ZZ/4 \oplus (\ZZ /2 )^{d'-1}\,\,{\rm if}\,\, k=0$$ Moreover
$$Pic_{tors} (X)\simeq (\QQ /\ZZ )^{d'-1}\oplus (\ZZ /2 )^{k'}$$

Assume $d$ is odd, $d:=2d'+1$. Thus $\bar{X}\setminus X$ is a real point, 
we have
$r=1$, $c=0$, $\eta (X)=0$. Moreover $g=d'$ and $k:=2k'+1$ is odd.
We have $t=k'$ and $s=k'+1$.
Then $$W(X)\simeq \ZZ^{k'+1}\oplus (\ZZ /2 )^{d'}$$ 
and
$$Pic_{tors} (X)\simeq (\QQ /\ZZ )^{d'}\oplus (\ZZ /2 )^{k'}$$

\subsubsection{Curves with equations $y^2 -P(x)=0$, $P$ monic}

Let $X$ be the plane curve with equation $y^2 -P(x)=0$. 
We denote by $d$ the degree of $P$ and $k$ the number of real roots 
of $P$.

First assume that $d$ is even, $d:=2d'$. We have $2$ real points
in $\bar{X}\setminus X$, hence $r=2$, $c=0$, $\eta (X)=0$ or $1$. 
We also get $g=d'-1$, $k=2k'$. If $k>0$, then 
$s=k'+1$, $t=k'-1$. If $k=0$ then $t=0$, $s=2$.
Then $$W(X)\simeq \ZZ^{k'+1}\oplus (\ZZ /2 )^{d'-1}\,\,{\rm if}\,\, k>0$$
$$W(X)\simeq \ZZ^{2}\oplus (\ZZ /2 )^{d'-1}\,\,{\rm if}\,\, k=0$$
Moreover
$$Pic_{tors} (X)\simeq (\QQ /\ZZ )^{d'-\eta(X)}\oplus (\ZZ /2 )^{k'-1}\,\,
{\rm if}\,\, k>0$$
and
$$Pic_{tors} (X)\simeq (\QQ /\ZZ )^{d'-\eta(X)}\,\,{\rm if}\,\, k=0$$

Assume $d:=2d'+1$ is odd. Then $r=1$, $c=0$, $\eta (X)=0$, $k=2k'+1$, $g=d'$,
$s=k'+1$, $t=k'$. We get
$$W(X)\simeq \ZZ^{k'+1}\oplus (\ZZ /2 )^{d'}$$
and
$$Pic_{tors} (X)\simeq (\QQ /\ZZ )^{d'}\oplus (\ZZ /2 )^{k'}$$

\subsubsection{Some remarks}

Let $X$ be the affine curve with equation $y^2 -P(x)$, with $P$ monic, 
and $X'$ the curve with equation $y^2 +P(x)$. 

If the degree $d$ of $P$ is odd, then we remark that we have obtained
the same results for $X$ and $X'$. This is not surprising:
by the isomorphism over $\RR$
$(x,y)\mapsto (-x,y)$, $X$ is isomorphic to the curve 
with equation $y^2 +Q(x)$ with $Q(x)=-P(-x)$ monic and 
$Q$ has exactly the same degree
and the same number of real roots than $P$.

In the remainder of this section we assume that the degree 
of $P$ is even, $d=2d'=2g+2$. We would like to know when $\eta(X)=0$ or $1$.
We have $\bar{X}\setminus X =\{P_1 ,P_2 \}=\bar{X_{\CC}}\setminus X_{\CC}$
and $\bar{X'}_{\CC}\setminus X'_{\CC}=\{ Q,\bar{Q}\}$ with $P_1 ,P_2$ 
real points 
and $Q$ a complex point.

\begin{prop}
Under the conditions stated above, we have $\eta(X)=\eta(X_{\CC} )=
\eta(X'_{\CC} )$.
\end{prop}

\begin{proof}
We have $\eta (X)=\eta (X_{\CC} )$ by Proposition \ref{comparaison}.
The curves $X'_{\CC}$ and $X_{\CC}$ are isomorphic over $\CC$
by $f:(x,y)\mapsto (x,iy)$. Hence
$Pic (X'_{\CC} )\simeq Pic( X'_{\CC} )$. Since they have the same number
of points at infinity, Proposition \ref{piccomplexe} shows that 
$\eta (X'_{\CC})=
\eta (X_{\CC} )$.
\end{proof}

Consequently $\eta (X)=1$ if and only if $P_1 -P_2$ is a torsion point
in the jacobian $Jac (\bar{X}) (\RR)=(Jac(\bar{X}_{\CC} )(\CC))^{G}$
if and only if $Q-\bar{Q}$ is a torsion point in
$Jac(\bar{X}_{\CC}) (\CC)$. More precisely, if these
points are torsion points in their Jacobian, they should have the same order
since the isomophism $f: X'_{\CC}\rightarrow X'_{\CC}; 
(x,y)\mapsto (x,iy)$ induces an isomorphism 
between $Jac (\bar{X'}_{\CC} )$ and $Jac(\bar{X}_{\CC})$ and 
the image of $Q-\bar{Q}$ is $P_1 -P_2$ interchanging $P_1$
and $P_2$ if necessary (\cite[Lem. 2.5]{Hu-Ma}).

We parametrize the set of genus $g$ hyperelliptic affine curves
like $X$ by $\RR^{2g+2}$ (the set of coefficients of $P$).
Let $M_{g,1}$ be the subset of curves with $\eta (X)=1$.
The following proposition may be proved in much the same ways as 
\cite[Lem. 2.4, Prop. 2.6]{Hu-Ma}.
\begin{prop}
The set $M_{g,1}$ has measure $0$ in $\RR^{2g+2}$.
\end{prop}

We will see that $M_{1,1} \not=\emptyset $ in the next section.

\subsubsection{Quartics}

Let $X$ be the smooth geometrically connected affine curve 
of equation $y^2 -P(x)=0$ with $P$ monic of degree $4$. 
We know that $g=1$ and we want to know when $\eta (X)=0$ or $1$.
Let $k=2k'$ be the number of real roots of $P$ and
we set $\bar{X}\setminus X=\{ P_1 ,P_2 \}$ as previously.
Then we distinguish $3$ cases: $k=0$, $k=2$, $k=4$.
In this section, we use ideas of \cite{Hu-Ma} that we adapt
to our problem.\\

{\bf Assume ${\bf k=0}$}. Then $X(\RR )$ has two connected components
and also $\bar{X} (\RR )$; moreover $P_1$ and $P_2$ are in two
different connected components of $\bar{X} (\RR )$. Therefore
$Jac (\bar{X}) (\RR )\simeq (\RR /\ZZ ) \oplus \ZZ /2 $
and $p=P_1 -P_2$ is in the non neutral component of 
$Jac (\bar{X}) (\RR )$ by \cite[Lem. 2.6]{S2}. Consequently
$\eta (X)=1$ if and 
only if $p$ is a torsion point of even order in $Jac (\bar{X}) (\RR )$.
We make an explicit calculation. Up to an isomorphism
we may assume that $P(x)=((x+b)^2+a^2 )((x+b)^2 +c^2 )$
with $a,b,c\in\RR$ and $a,c>0$. We get 
a Weierstrass equation 
$u^2 =(v +4b^2 )(u -(c-a)^2 )(v -(c+a)^2 )$
Then $p$ has coordinates $(0, 2b(c^2 -a^2))$ and we have 
$3$ points of order $2$: $p_1 =(-4b^2 ,0 )$, $p_2 =((c-a)^2 ,0))$,
$p_3 =((c+a)^2 ,0))$. 
We see that $p,p_1 ,p_2$ are on the non neutral component.
Thus
\begin{prop} Under the above conditions:
$$Pic_{tors} (X)\simeq (\QQ /\ZZ )$$
if and only if $(2n-1)p=p_1$ or $(2n-1)p=p_2$ or $2np =p_3$
for a $n\in\NN\setminus \{ 0\}$.
Else 
$$Pic_{tors} (X)\simeq (\QQ /\ZZ )^2$$
\end{prop}
When $\eta(X)=1$, using the duplication formula,
each case $(2n-1)p=p_1$ or $(2n-1)p=p_2$ or $2np =p_3$ 
is equivalent to a polynomial equation 
in $a,b,c$ defined over $\QQ$.
For example, $p=p_1$ if and only if $b=0$,
$p=p_2$ if and only if $c=a$ and 
$2p=p_3$ if and only if $-16b^2 ac+c^4 -2a^2 c^2 +a^4 =0$.\\

{\bf Assume ${\bf k=2}$}. Then $X(\RR )$ has two connected components
but $\bar{X} (\RR )$ has only one. We get
$Jac (\bar{X}) (\RR )\simeq \RR /\ZZ $
and $p=P_1 -P_2$ is in the neutral component of 
$Jac (\bar{X}) (\RR )$.
Up to an isomorphism
we may assume that $P(x)=((x+b)^2+a^2 )((x-b)^2 -c^2 )$
with $a,b,c\in\RR$ and $a,c>0$. We get 
a Weierstrass equation 
$u^2 =(v +4b^2 )(v^2 -2(c^2 -a^2)v +(c^2+a^2)^2 )$.
Then $p$ has coordinates $(0, 2b(c^2 +a^2))$ and we have only
$1$ point of order $2$: $p_1 =(-4b^2 ,0 )$. 
\begin{prop} Under the above conditions:
$$Pic_{tors} (X)\simeq (\QQ /\ZZ )$$
if and only if $np=p_1$ or $2np=-p$
for a $n\in\NN\setminus \{ 0\}$.
Else 
$$Pic_{tors} (X)\simeq (\QQ /\ZZ )^2$$
\end{prop}
For example $p=p_1$ if and only if $b=0$,
$2p=p_1$ if and only if 
$$-64c^4 b^4+128 a^2 b^4 c^2-64 a^4 b^4-a^8-c^8+1
6 b^2 c^6-16 b^2 a^6$$ $$-4 a^6 c^2-6 a^4 c^4-4 a^2 c^6+
16 b^2 c^4 a^2-16 b^2 a^4 c^2 =0$$
and $2p=-p$ if and only if 
$$a^8+c^8-16 b^2 c^6+16 b^2 a^6+
4 a^6 c^2+6 a^4 c^4+4 a^2 c^6-16 b^2 c^4 a^2+
16 b^2 a^4 c^2-256 a^2 b^4 c^2 =0$$\\

{\bf Assume ${\bf k=4}$}. Then $X(\RR )$ has 3 connected components
but $\bar{X} (\RR )$ has only 2. We get
$Jac (\bar{X}) (\RR )\simeq (\RR /\ZZ )\oplus \ZZ /2$
and $p=P_1 -P_2$ is in the neutral component of 
$Jac (\bar{X}) (\RR )$ since $P_1$ and $P_2$ lie
to the same connected component of $\Bar{X} (\RR)$.
Up to an isomorphism
we may assume that $P(x)=((x+b)^2 -a^2 )((x-b)^2 -c^2 )$
with $a,b,c\in\RR$ and $a,c>0$. We get 
a Weierstrass equation 
$u^2 =(v +4b^2 )(v +(c-a)^2 )(v +(c+a)^2 )$
Then $p$ has coordinates $(0, 2b(c^2 -a^2))$ and we have 
$3$ points of order $2$: $p_1 =(-4b^2 ,0 )$, $p_2 =(-(c+a)^2 ,0))$,
$p_3 =(-(c-a)^2 ,0))$. We see that $p,p_3$ 
are on the neutral component if $ 4b^2 > (c-a)^2 $
and $p,p_1$ 
are on the neutral component if $ 4b^2 < (c-a)^2 $.
We have to remark that $2b\not= \pm (c-a)$ and 
$2b\not= \pm (c+a)$ since the discriminant of $P$
is $\not= 0$ (the curve is smooth).
\begin{prop} Under the above conditions:
$$Pic_{tors} (X)\simeq (\QQ /\ZZ )\oplus (\ZZ /2 )$$
if and only $np=p_3$ or $2np=-p$ in the case $ 4b^2 > (c-a)^2 $,
and $np=p_1$ or $2np=-p$ in the case $ 4b^2 < (c-a)^2 $,
for a $n\in\NN\setminus \{ 0\}$.
Else 
$$Pic_{tors} (X)\simeq (\QQ /\ZZ )^{2}\oplus (\ZZ /2 )$$
\end{prop}
Assume $ 4b^2 > (c-a)^2 $, then for example, 
$p=p_3$ if and only if $a=c$,
$2p=-p$ if and only if
$$16c^4 b^2 a^2+16 c^2 a^4 b^2+256 b^4 c^2 a^2+
c^8+a^8$$ $$-4 c^6 a^2+6 c^4 a^4-16 c^6 b^2-4 c^2 a^6-16 a^6 b^2 =0$$
Assume $ 4b^2 > (c-a)^2 $ then $p=p_1$ if and only if
$b=0$.\\
\\

Following \cite{Hu-Ma}, if we assume that $P$ admits
a factorization $P(x)=((x+b)^2 +a^2 )((x-b)^2 +c^2 )$
in the case $k=0$, $P(x)=((x+b)^2 +a^2 )((x-b)^2 -c^2 )$
in the case $k=2$, $P(x)=((x+b)^2 -a^2 )((x-b)^2 -c^2 )$
in the case $k=4$; with $a,b,c\in\QQ$ and $a,c>0$.
Then $p$ and all the points denoted by $p_1 ,p_2 ,p_3$
are rational points. Now we used a famous theorem of Mazur
wich asserts that the torsion subgroup of $Jac (\bar{X}) (\QQ)$ 
is either $\ZZ / n\ZZ $ for $n=1,2,3,\ldots,10,12$
or $\ZZ /n\ZZ \oplus\ZZ /2\ZZ$ for $n=2,4,6,8$.
We see that if $k=0$ or $k=4$, $(Jac (\bar{X}) (\QQ))_{tors}$ is
of the second type since it has already $3$ distinct elements
of order $2$. If $k=2$, $(Jac (\bar{X}) (\QQ))_{tors}$ is clearly of the first 
type with $n$ even.
Therefore we obtain a finite number of conditions
to assert that $\eta (X)=1$.
\begin{prop}
We assume that $P$ admit a rational factorization. Then,
\begin{description}
\item[({\cal{i}})] 
If $k=0$ then $\eta (X)=1$ if and only if
$(2n-1)p=p_1$ or $(2n-1)p=p_2$ or $2np =p_3$
for  $0<n\leq 2$ ($6$ cases).
\item[({\cal{ii}})] If $k=2$ then $\eta (X)=1$ if and only if
$np=p_1$ for $n=1,\ldots,6$ or $2np=-p$
for $n=1,2,3,4$ ($10$ cases).
\item[({\cal{iii}})] If $k=4$ and $ 4b^2 > (c-a)^2 $ 
then $\eta (X)=1$ if and only if
$np=p_3$ for $n=1,\ldots,4$ ($4$ cases) ($2np=-p$
is not allowed).
\item[({\cal{ii}})] If $k=4$ and $ 4b^2 < (c-a)^2 $ 
then $\eta (X)=1$ if and only if
$np=p_1$ for $n=1,\ldots,4$ ($4$ cases).
\end{description}
\end{prop}

\end{document}